\documentclass[11pt]{amsart}
\usepackage{amsmath}
\usepackage{amssymb}
\usepackage{amsthm}
\usepackage{array}
\usepackage{xy}
\usepackage[pdftex]{graphicx}
\usepackage{color}
\usepackage{transparent}
\usepackage{latexsym}
\usepackage{hyperref}

\newtheorem{theorem}{Theorem}[section]
\newtheorem{proposition}[theorem]{Proposition}
\newtheorem{lemma}[theorem]{Lemma}
\newtheorem{corollary}[theorem]{Corollary}
\theoremstyle{remark}

\numberwithin{equation}{section}

\DeclareMathOperator{\diam}{diam}
\DeclareMathOperator{\dist}{dist}

\DeclareMathOperator{\spt}{spt}
\DeclareMathOperator{\aff}{aff}
\DeclareMathOperator{\ri}{ri}

\newcommand{\eps}{\varepsilon}
\newcommand{\dd}{\,\mathrm d}
\allowdisplaybreaks[2]
\hypersetup{hidelinks,pdftitle={Global regularity for optimal transportation on c-convex domains},pdfauthor={Yuanyuan Li, Xianduo Wang and Xu-Jia Wang}}

\title{Global regularity in optimal transportation\\under the A3w condition}

\author{Yuanyuan Li}
\address{Institute for Theoretical Sciences, Westlake University, Hangzhou 310030, China}
\email{lyyuan@westlake.edu.cn}
\author{Xianduo Wang}
\address{School of Mathematics and Statistics, Huazhong University of Science and Technology, Wuhan 430074, China}
\email{wangxd@hust.edu.cn}
\author{Xu-Jia Wang}
\address{Institute for Theoretical Sciences, Westlake University, Hangzhou 310030, China}
\email{wangxujia@westlake.edu.cn}
\date{}

\begin{document}
\begin{abstract}
We prove global regularity estimates for optimal transportation under A3w on bounded mutually $c$-convex domains. For nondegenerate continuous densities, we prove global $W^{2,p}$ regularity of the
potentials up to the boundary. If in addition the boundary is $C^{1,\beta}$ and the densities are $C^{\alpha}$, we improve this to $C^{2,\gamma}$, where $\gamma=\min\{\alpha,\beta\}$. 
\end{abstract}
\maketitle

\section{Introduction}

Let $\Omega,\Omega^*\subset\mathbb R^n$, $n\ge3$, be two bounded domains with densities $f$ and $g$ satisfying $\int_\Omega f\,dx=\int_{\Omega^*}g\,dy.$ Assume that for some $0<\lambda<1$,
\begin{equation}\label{density}
\begin{gathered}
\lambda\le f\le\lambda^{-1}\quad\text{on }\overline\Omega,\qquad \lambda\le g\le\lambda^{-1}\quad\text{on }\overline{\Omega^*}.
\end{gathered}
\end{equation}
The optimal transportation problem is to maximize the cost functional
\[
\mathcal C(T)=\int_\Omega c(x,T(x))f(x)\,dx
\]
among maps $T:\Omega\to\Omega^*$ satisfying $T_\#(f\,dx)=g\,dy$.

Throughout, we assume $c\in C^4(\mathcal U)$ for some neighborhood $\mathcal U$ of $\overline\Omega\times\overline{\Omega^*}$. We also assume the following twist and nondegeneracy conditions.
\begin{itemize}
\item[\textup{(A1)}] For each $x\in\overline\Omega$, the mapping $y\mapsto c_x(x,y)$ is one-to-one on $\overline{\Omega^*}$. For each $y\in\overline{\Omega^*}$, the mapping $x\mapsto c_y(x,y)$ is one-to-one on $\overline\Omega$.
\item[\textup{(A2)}] The matrix $c_{xy}$ is invertible on $\overline\Omega\times\overline{\Omega^*}$.
\end{itemize}
Under \textup{(A1)}, there is an optimal mapping $T$ satisfying $c_x(x,T(x))=Du(x)$, where $u$ is its $c$-convex potential. We denote the dual potential by $v$. See Gangbo~\cite{Gangbo1995} and De Philippis and Figalli~\cite{DPF2015}.

Let $Y(x,\cdot)$ denote the inverse of $y\mapsto c_x(x,y)$ and set $A(x,p)=c_{xx}(x,Y(x,p))$. For a smooth optimal diffeomorphism $T(x)=Y(x,Du(x))$, one can get the Monge--Amp\`ere type equations with natural boundary condition
\begin{equation}\label{ma}
\det[D^2u-A(x,Du)]=|\det c_{xy}(x,T(x))|\frac{f(x)}{g(T(x))},\quad
T(\Omega)=\Omega^*.
\end{equation}

For $c(x,y)=x\cdot y$, by Brenier's theorem~\cite{Brenier1991}, $T=Du$ with $u$ convex, and \eqref{ma} reduces to the Monge--Amp\`ere equation. Caffarelli~\cite{Caffarelli1990W2p,Caffarelli1991,Caffarelli1992Mappings} developed the interior theory. For interior Sobolev estimates with right hand side bounded from above and below, see~\cite{PhillipisFigalli2013} and~\cite{PhillipisFigalliSavin2013}. For the classical boundary estimates, see Delano\"e, Caffarelli and Urbas~\cite{Dela,Caffarelli1992Boundary,Caffarelli1996Boundary,Urbas}. Chen, Liu and Wang~\cite{CLW2021} proved global $W^{2,p}$ estimates for continuous densities and global $C^{2,\alpha}$ estimates for $C^\alpha$ densities, assuming the convex domains with $C^{1,1}$ boundaries, 

For arbitrary bounded convex domains, the global Sobolev estimate was obtained by Savin and Yu~\cite{SavinYu2020} for constant densities in the plane, and in all dimensions by Collins and Tong~\cite[Theorem 1.1]{CT2026} for Dini continuous densities and by Chen, Li and Wang~\cite[Theorem 1.1]{CLW2026} for continuous densities. For $C^{1,\beta}$ convex domains and $C^\alpha$ densities, Collins and Tong~\cite[Theorem 1.2]{CT2026} obtained global $C^{2,\min\{\alpha,\beta\}}$ regularity.

For general cost functions, Ma, Trudinger and Wang~\cite{MTW2005} proved interior $C^3$ regularity under the \textup{(A3s)} for positive $C^2$ densities and a $c^*$-convex target domain, where \textup{(A3s)} requires a uniform positive lower bound in the inequality below for (A3w). 
\begin{itemize}
\item[\textup{(A3w)}] For $x\in\overline\Omega$, $p\in c_x(x,\overline{\Omega^*})$ and vectors $\xi,\eta\in\mathbb R^n$ with $\xi\cdot\eta=0$,
\[
D^2_{p_kp_l}A_{ij}(x,p)\,\xi_i\xi_j\eta_k\eta_l\ge0.
\]
\item[\textup{(B)}] The domain $\Omega$ is $c$-convex with respect to $\overline{\Omega^*}$, and $\Omega^*$ is $c^*$-convex with respect to $\overline\Omega$.
\end{itemize}

Loeper~\cite{Loeper2009} proved that failure of \textup{(A3w)} can lead to discontinuity of the optimal mapping even for smooth positive densities. Kim and McCann~\cite{KM2010} interpreted \textup{(A3w)} in terms of the sectional curvature of a pseudo-Riemannian metric. Under \textup{(A3w)} and the uniform domain convexity assumptions of Figalli, Kim and McCann~\cite[Theorem 2.1]{FKM2013}, interior $C^{2,\alpha}$ regularity for H\"older densities follows from their strict $c$-convexity and differentiability theorem combined with the interior estimate of Liu, Trudinger and Wang~\cite[Theorem 1.2 and Remark 4.1]{LTW2010}.

For global regularity, Trudinger and Wang~\cite{TW2009} proved global $C^3$ regularity for \eqref{ma} under \textup{(A3w)}, with additional hypotheses for the second derivative estimate. Trudinger~\cite{Trudinger2013} later removed the additional conditions for positive $C^2$ densities on mutually uniformly $c$-convex $C^4$ domains. Hong and Wang~\cite{HongWang2026} proved global $C^{2,\min\{\alpha,\beta\}}$ regularity for smooth costs satisfying \textup{(A3s)}, positive $C^\alpha$ densities and mutually $c$-convex domains with $C^{1,\beta}$ boundaries. We assume one of the following conditions on the densities.
\begin{itemize}
\item[\textup{(C)}] $f\in C(\overline\Omega)$ and $g\in C(\overline{\Omega^*})$.
\item[\textup{(H)}] $f\in C^{\alpha}(\overline\Omega)$ and $g\in C^{\alpha}(\overline{\Omega^*})$ for some $0<\alpha<1$.
\end{itemize}

Our main results are as follows. 

\begin{theorem}\label{SR}
Assume \textup{(A1)}, \textup{(A2)}, \textup{(A3w)}, \textup{(B)}, \textup{(C)} and~\eqref{density}. After fixing the additive constants, then
\begin{equation}\label{sobolevbound}
\|u\|_{W^{2,p}(\Omega)}+\|v\|_{W^{2,p}(\Omega^*)}\le C_p,
\qquad 1<p<\infty.
\end{equation}
The map $T$ extends to a homeomorphism of $\overline\Omega$ onto $\overline{\Omega^*}$.
\end{theorem}

\begin{theorem}\label{HR}
Under the hypotheses of Theorem~\ref{SR}, assume \textup{(H)} 
\[
\partial\Omega,\partial\Omega^*\in C^{1,\beta},
\qquad 0<\beta<1.
\]
After fixing the additive constants, then we have
\begin{equation}\label{holderbound}
\|u\|_{C^{2,\gamma}(\overline\Omega)}+\|v\|_{C^{2,\gamma}(\overline{\Omega^*})}\le C,
\end{equation}
where $\gamma=\min\{\alpha,\beta\}$. The map $T$ extends to a $C^{1,\gamma}$ diffeomorphism of $\overline\Omega$ onto $\overline{\Omega^*}$.
\end{theorem}

Under the hypotheses of Theorem~\ref{SR}, if $\partial\Omega,\partial\Omega^*\in C^{1,\beta}$ for some $0<\beta<1$, Proposition~\ref{oblique} gives
\begin{equation}\label{angle}
\nu_\Omega(x)^Tc_{xy}(x,T(x))^{-T}\nu_{\Omega^*}(T(x))\ge C^{-1}>0,
\qquad x\in\partial\Omega,
\end{equation}
where $\nu_\Omega$ and $\nu_{\Omega^*}$ are the inward unit normals.

The main step is to control boundary sections without assuming obliqueness or differentiability of the boundaries. For general cost functions, two difficulties will arise. The radial functions $x\cdot Du(x)$ and $y\cdot Dv(y)$ need not agree on the transport graph. Moreover, affine normalization may enlarge the higher order remainder of the cost when the sections are highly eccentric.

We compare sections of the potentials with the same mass. In the coordinates~\eqref{cnorm}, let $h(m)$ and $h^*(m)$ be determined by
\[
\int_{S_{h(m)}}f\,dx=\int_{S^*_{h^*(m)}}g\,dy=m.
\]
Observe that section $S^*_{h^*(m)}$ minimizes the integral of $v$ among all target sets of mass $m$, with minimum $\int_0^m h^*(s)\,ds$. Mass conservation therefore gives a comparison with the image of $S_{h(m)}$. This is the reason why we introduce the inverse distribution functions of the section masses.

By \textup{(A3w)}, we obtain a comparison for $c$-functions along rays. Combining radial integration with the rearrangement inequality on both sides, we obtain
\[
\begin{gathered}
J(m)=m^{-1-2/n}\int_0^m\bigl(h(s)+h^*(s)\bigr)\,ds,
\qquad J(m)\ge C^{-1},\\
\frac{d\log J(m)}{d\log m}\ge-\omega(m),
\qquad 0\le\omega(m)\longrightarrow0
\quad \text{ as } m\rightarrow0.
\end{gathered}
\]
The normalization of $J$ is chosen so that it is constant for potentials homogeneous of degree two with constant densities. Note that the radial integration uses only first derivatives of the potentials and we do not need approximation argument to insure the regularity used in calculation.

For continuous densities, we use the scale counting argument used by Chen, Li and Wang~\cite[Propositions 7.2 and 7.3]{CLW2026}. By iteration, we obtain
\[
c_\varepsilon h^{1/2+\varepsilon}
\le \sigma_{\min}(L_h)
\le \|L_h\|
\le C_\varepsilon h^{1/2-\varepsilon}.
\] 
The same estimate holds for $L_h^*$. With $P_h=hL_h^{-T}$, we have
\[
\bigl\|h^{-1}c(L_hX,P_hY)-X\cdot Y\bigr\|_{C^{2,1}(B_R\times B_R)}
\longrightarrow 0
\]
for every fixed $R<\infty$. Applying the interior estimate of Chen and Figalli~\cite{CF2017}, we get Theorem~\ref{SR}. 

For $C^{1,\beta}$ boundaries, the growth estimates imply uniform obliqueness. Under \textup{(H)}, conducting a convex correction and reflection, allow comparison with solutions of the Monge--Amp\`ere equation with constant right hand side. The first iteration gives positive definite quadratic expansions at boundary points. These expansions provide the quadratic scale needed for a second iteration, which gives the exponent $\gamma=\min\{\alpha,\beta\}$. Combining the boundary estimates with the interior estimate of De Philippis and Figalli~\cite[Theorem 5.3]{DPF2015} proves Theorem~\ref{HR}.

\section{Preliminaries}\label{prelim}

Unless otherwise stated, we assume \textup{(A1)}, \textup{(A2)}, \textup{(A3w)}, \textup{(B)}, \textup{(C)} and~\eqref{density} throughout the paper.

For $x_0\in\overline\Omega$ and $y_0\in\overline{\Omega^*}$, we say that $x\mapsto u(x_0)+c(x,y_0)-c(x_0,y_0)$ is a $c$-support of $u$ at $x_0$ if it lies below $u$ on $\overline\Omega$. The $c$-normal mapping $T_u$ at $x_0$ is the set of all such $y_0$~\cite[Section 2.1]{LTW2015}. We say that $u$ is $c$-convex if $T_u(x)$ is nonempty for every $x\in\overline\Omega$, and strictly $c$-convex if, in addition, each $c$-support touches $u$ at only one point. For any Borel set $E\subset\overline\Omega$, we write $T_u(E)=\bigcup_{x\in E}T_u(x)$. We use the same definitions for $v$ with the dual cost $c^*(y,x)=c(x,y)$ and denote its $c^*$-normal mapping by $T_v$.

The two potentials are dual to each other and satisfy
\begin{equation}\label{duality}
\begin{aligned}
u(x)&=\max_{y\in\overline{\Omega^*}}\{c(x,y)-v(y)\},\\
v(y)&=\max_{x\in\overline\Omega}\{c(x,y)-u(x)\}.
\end{aligned}
\end{equation}
Thus $u(x)+v(y)\ge c(x,y)$, with equality if and only if $y\in T_u(x)$. We call such a pair $(x,y)$ a contact pair. Let
\[
\Gamma=\{(x,y)\in\overline\Omega\times\overline{\Omega^*}\mid u(x)+v(y)=c(x,y)\}
\]
be the closed contact set. The optimal plan $\pi=(\operatorname{Id},T)_\#(f\,dx)$ satisfies
\[
\spt\pi\subset\Gamma,\qquad \pi_x=f\,dx,\qquad \pi_y=g\,dy.
\]
Before uniqueness of contacts is proved, $T_u$ always denotes the entire set-valued $c$-normal mapping.

For fixed $y_0\in\overline{\Omega^*}$, a $c$-segment with respect to $y_0$ is a curve whose image under $x\mapsto c_y(x,y_0)$ is a straight segment. The domain $\Omega$ is $c$-convex with respect to $\overline{\Omega^*}$ if $c_y(\Omega,y)$ is convex for every $y\in\overline{\Omega^*}$. For $C^2$ boundaries, uniform $c$-convexity means that the principal curvatures of the boundaries of these convex domains have a uniform positive lower bound. The corresponding definitions for $\Omega^*$ use the dual cost $c^*$.

In the derivatives of $c$, subscripts before the comma refer to $x$ and those after the comma refer to $y$. Let $c^{r,k}$ denote the $(r,k)$ entry of $(c_{xy})^{-1}$, so that $c_{i,r}c^{r,k}=\delta_i^k$. As in~\cite[(1.19)]{TW2009}, condition \textup{(A3w)} can be written at $(x,Y(x,p))$ as
\begin{equation}\label{mtw}
(c_{ij,rs}-c_{ij,a}c^{a,b}c_{b,rs})c^{r,k}c^{s,l}\xi_i\xi_j\eta_k\eta_l\ge0,
\qquad \xi\cdot\eta=0.
\end{equation}
Repeated indices are summed from $1$ to $n$.

We write
\[
\omega_f(t)=\sup\{|f(x)-f(x')|\mid x,x'\in\overline\Omega,\ |x-x'|\le t\}
\]
and define $\omega_g$ similarly on $\overline{\Omega^*}$. Unless otherwise stated, constants depend only on $n,\lambda$, the $C^4$ norm of $c$, the bound for $c_{xy}^{-1}$, the domains, and the moduli $\omega_f,\omega_g$. In Theorem~\ref{HR}, they may also depend on $\alpha,\beta$, the density norms, and the $C^{1,\beta}$ boundary charts. The letter $C$ may change from line to line. We write $C_p$ when the constant also depends on $p$.

The Aleksandrov identity in \cite[Definition 2.1 and (2.2)]{LTW2015} gives
\begin{equation}\label{mass}
\int_{T_u(E)\cap\Omega^*}g\dd y=\int_{E\cap\Omega}f\dd x
\qquad\text{for every Borel }E\subset\overline\Omega.
\end{equation}
Here images are measured in the completed Lebesgue measure.
Indeed, $v$ is locally semiconvex and hence differentiable almost everywhere in $\Omega^*$.
At such a point,
\[
(x,y)\in\Gamma\quad\Longrightarrow\quad Dv(y)=c_y(x,y),
\]
so the twist assumption makes the inverse contact unique.
Disintegration of $\pi$ therefore gives \eqref{mass}, including sets meeting $\partial\Omega$.
The dual identity also holds.

Fix $(x_0,y_0)\in\Gamma$ and write the original variables as $\bar x,\bar y$.
Set $M=c_{xy}(x_0,y_0)$ and introduce the $c$-exponential coordinates
\[
x=c_y(\bar x,y_0)-c_y(x_0,y_0),
\qquad
y=M^{-1}\bigl(c_x(x_0,\bar y)-c_x(x_0,y_0)\bigr).
\]
Let $\bar x=\bar x(x)$ and $\bar y=\bar y(y)$ be their inverses.
As in \cite[Section 2.2, (2.4)--(2.8)]{LTW2015}, put
\[
\begin{aligned}
\tilde c(x,y)&=c(\bar x,\bar y)-c(\bar x,y_0)-c(x_0,\bar y)+c(x_0,y_0),\\
\tilde u(x)&=u(\bar x)-u(x_0)-c(\bar x,y_0)+c(x_0,y_0),\\
\tilde v(y)&=v(\bar y)-v(y_0)-c(x_0,\bar y)+c(x_0,y_0).
\end{aligned}
\]
The densities in these coordinates are
\[
\tilde f(x)=f(\bar x(x))|\det D_x\bar x|,
\qquad
\tilde g(y)=g(\bar y(y))|\det D_y\bar y|.
\]
Both coordinate domains are convex.
We have
\begin{equation}\label{cnorm}
\begin{gathered}
\tilde c(x,0)=\tilde c(0,y)=0,
\qquad \tilde c_y(x,0)=x,
\qquad \tilde c_x(0,y)=y,\\
\tilde u(0)=\tilde v(0)=0,
\qquad \tilde u,\tilde v\ge0.
\end{gathered}
\end{equation}
For the rest of this section, we retain $c,u,v,f,g,\Omega,\Omega^*$ for the transformed quantities.
The coordinate changes are $C^3$ with uniformly bounded derivatives and inverse derivatives.
The transformed cost has bounded third derivatives and bounded mixed derivatives $D_x^2D_y^2c$.
Condition A3w, mutual $c$-convexity and the mass identity are preserved.
These bounds are uniform in $(x_0,y_0)\in\Gamma$.

Define the sections
\[
S_h=\{x\in\Omega\mid u(x)<h\},
\qquad
S_h^*=\{y\in\Omega^*\mid v(y)<h\}.
\]
They are convex by Loeper's maximum principle \cite[Theorem 3.2]{Loeper2009}, in the form proved under the present domain assumptions by Kim and McCann \cite[Theorem 4.10]{KM2010}.
Mutual $c$-convexity supplies the segments required in that theorem.
The boundary cases follow by continuity.
John's ellipsoid theorem \cite{John1948}, in the form \cite[(2.9)]{LTW2015}, gives positive definite matrices $L_h,L_h^*$ such that
\[
\begin{gathered}
x_h+L_hB_1\subset S_h\subset x_h+nL_hB_1,\\
y_h+L_h^*B_1\subset S_h^*\subset y_h+nL_h^*B_1.
\end{gathered}
\]
Since $0\in\overline{S_h}\cap\overline{S_h^*}$,
\[
|L_h^{-1}x_h|\le n,
\qquad |(L_h^*)^{-1}y_h|\le n.
\]

By \eqref{cnorm} and \cite[(2.8)]{LTW2015},
\[
A(x,0)=0,\qquad D_pA(x,0)=0,\qquad |D^2_{pp}A|\le C.
\]
For $p\in c_x(x,\overline{\Omega^*})$, condition \textup{(B)} gives $[0,p]\subset c_x(x,\overline{\Omega^*})$. Hence
\[
A(x,p)[\xi,\xi]
=\int_0^1(1-t)D^2_{pp}A(x,tp)[\xi,\xi,p,p]\,\dd t.
\]
For $p\ne0$, let
\[
\xi^\perp=\xi-\frac{\xi\cdot p}{|p|^2}p,
\qquad \xi^\perp\cdot p=0,
\qquad |\xi^\perp|\le|\xi|.
\]
By \textup{(A3w)},
\[
A(x,p)[\xi^\perp,\xi^\perp]\ge0,
\qquad |A(x,p)|\le C|p|^2.
\]
Therefore
\begin{equation}\label{scalar}
\begin{aligned}
A(x,p)[\xi,\xi]
&\ge-C\bigl(|\xi^\perp||p||\xi\cdot p|+|\xi\cdot p|^2\bigr)\\
&\ge-C|\xi||p||\xi\cdot p|.
\end{aligned}
\end{equation}
For a $c$-function $\varphi$ with focus $y\in\overline{\Omega^*}$,
\[
D^2\varphi[\xi,\xi]\ge-C|\xi||D\varphi||D\varphi\cdot\xi|,
\qquad |D\varphi|\le C|y|.
\]

For $h>0$ and an invertible matrix $L$, set
\[
w(z)=h^{-1}\varphi(x_b+Lz).
\]
Then
\begin{equation}\label{rescaledscalar}
\begin{aligned}
D^2w[e,e]
&=h^{-1}D^2\varphi[Le,Le]\\
&\ge-C|Le||D\varphi||Dw\cdot e|
\ge-C\|L\||e||Dw\cdot e|.
\end{aligned}
\end{equation}
Here $C$ is independent of $h$ and $L$.

The one-dimensional inequality $\psi''\ge-\kappa|\psi'|$ gives
\[
\psi'(s)>0
\quad\Longrightarrow\quad
\psi'(t)\ge e^{-\kappa(t-s)}\psi'(s)
\qquad(t\ge s).
\]
For fixed $r,R>0$ and bounded $\|L\|$, suppose
\[
x_b+LB_r(z_0)\subset\Omega,\qquad
\sup_{B_r(z_0)}w\le H,\qquad w(z_0)\ge-H.
\]
Then
\begin{equation}\label{propagation}
|Dw(z_0)|\le CH,\qquad w(z)\ge-CH,
\end{equation}
whenever
\[
|z-z_0|\le R,\qquad
[x_b+Lz_0,x_b+Lz]\subset\overline\Omega.
\]
A $c$-function with an interior maximum on a segment is constant on that segment. Applying the same comparison to a $c$-support of $u$ at $tx$, we obtain, for a uniform $b\ge1$,
\begin{equation}\label{radial}
\begin{gathered}
u(tx)\le bt\,u(x)
\qquad(x\in\overline\Omega,\ 0\le t\le1),\\
S_{ah}\subset baS_h,\qquad
|S_{ah}|\le(ba)^n|S_h|
\qquad(a\ge1).
\end{gathered}
\end{equation}
The same estimates hold for the dual problem.

\begin{proposition}\label{volume}
There exist $C>1$ and $h_0>0$, independent of the contact pair, such that, for $0<h<h_0$,
\[
C^{-1}|S_h|\le|S_h^*|\le C|S_h|,
\]
\[
C^{-1}h^n\le|S_h|\le Ch^{n/2},
\qquad
C^{-1}h^n\le|S_h^*|\le Ch^{n/2}.
\]
\end{proposition}

\begin{proof}
Let $E_h=x_h+L_hB_{1/4}$. For $y\in T_u(E_h)$, choose $z_0\in B_{1/4}$ such that $(x_h+L_hz_0,y)\in\Gamma$, and set
\[
w(z)=h^{-1}\bigl(c(x_h+L_hz,y)-v(y)\bigr).
\]
Then
\[
w\le1\quad\text{in }B_1,\qquad
w(z_0)=h^{-1}u(x_h+L_hz_0)\ge0.
\]
Since $B_{1/2}(z_0)\subset B_1$, we may apply \eqref{propagation} at $z_0$ and then at $0$ to obtain
\[
w(0)\ge-C,\qquad |Dw(0)|\le C,
\qquad
w(-L_h^{-1}x_h)=-\frac{v(y)}h\ge-C.
\]
Consequently,
\[
T_u(E_h)\cap\Omega^*\subset S_{Ch}^*,
\qquad
|L_h^Tc_x(x_h,y)|\le Ch
\quad(y\in T_u(E_h)).
\]
By \eqref{mass}, the density bounds and \eqref{radial},
\[
\begin{aligned}
|S_h|
&\le C|E_h|
\le C|T_u(E_h)\cap\Omega^*|\\
&\le C|S_{Ch}^*|
\le C|S_h^*|.
\end{aligned}
\]

By \textup{(A1)}--\textup{(A2)}, the change of variables $y\mapsto L_h^Tc_x(x_h,y)$ gives
\[
C^{-1}\det L_h
\le |T_u(E_h)\cap\Omega^*|
\le C\frac{h^n}{\det L_h}.
\]
Thus
\[
(\det L_h)^2\le Ch^n,\qquad
|S_h|\le C\det L_h\le Ch^{n/2}.
\]

We also have
\[
\|Du\|_{L^\infty(\Omega)}
+\|Dv\|_{L^\infty(\Omega^*)}\le C.
\]
The coordinate domains are convex with uniform inradius and diameter bounds. Since $u(0)=0$,
\[
\Omega\cap B_{h/C}\subset S_h,\qquad
|S_h|\ge|\Omega\cap B_{h/C}|\ge C^{-1}h^n.
\]
Interchanging $u$ and $v$ completes the proof.
\end{proof}

\begin{proposition}\label{contacts}
The $c$-normal mappings $T_u$ and $T_v$ are single-valued on the closures and are inverse homeomorphisms between $\overline\Omega$ and $\overline{\Omega^*}$. Moreover,
\[
u\in C^1(\Omega),\qquad v\in C^1(\Omega^*),
\]
\[
Du(x)=c_x(x,T_u(x)),\qquad
Dv(y)=c_y(T_v(y),y),
\]
with
\[
Du\in C(\overline\Omega;\mathbb R^n),\qquad
Dv\in C(\overline{\Omega^*};\mathbb R^n).
\]
The sections satisfy
\[
\lim_{h\rightarrow0}\sup_{(x_0,y_0)\in\Gamma}
\bigl(\diam S_h+\diam S_h^*\bigr)=0.
\]
\end{proposition}

\begin{proof}
Suppose that $T_v(y_0)$ contains more than one point. Use the coordinates \eqref{cnorm} with focus $y_0$ and set
\[
E=\{x\in\overline\Omega\mid u(x)=0\}.
\]
The set $E$ is convex and $\diam E>0$. Choosing the source origin in $\ri E$ and making the corresponding change of target coordinates, we may assume
\[
0\in\ri E.
\]
Let
\[
F=\{y\in\overline{\Omega^*}\mid v(y)=0\}.
\]
Then $F$ is convex and $0\in F$. By duality,
\[
c(z,y)\le u(z)+v(y)=0=c(0,y)
\qquad(z\in E,\ y\in F).
\]
Since $0\in\ri E$, the one-dimensional comparison gives
\[
c=0\quad\text{on }E\times F.
\]

For $y\in\ri F$ and $(x,y)\in\Gamma$,
\[
c(x,z)-u(x)\le v(z)=0
\quad(z\in F),\qquad
c(x,y)-u(x)=0.
\]
The same comparison on $F$ gives
\[
c(x,z)-u(x)=0\qquad(z\in F).
\]
Taking $z=0$, we obtain
\begin{equation}\label{inversecontact}
y\in\ri F,\quad(x,y)\in\Gamma
\quad\Longrightarrow\quad x\in E.
\end{equation}

Let $V=\aff E$, and let $\operatorname{proj}_V$ denote orthogonal projection onto $V$. Since $0\in\ri E$ and $c_x(0,y)=y$, the one-dimensional comparison gives, for $y\in\overline{\Omega^*}$,
\[
\max_{z\in E}c(z,y)\ge C^{-1}|\operatorname{proj}_V y|,
\qquad
\min_{z\in E}c(z,y)\ge-C|\operatorname{proj}_V y|.
\]
Here $C$ may also depend on $E$. For $y\in S_h^*$,
\[
c(z,y)\le u(z)+v(y)=v(y)<h
\qquad(z\in E).
\]
Hence
\begin{equation}\label{strip}
|\operatorname{proj}_V y|\le Ch,\qquad
-Ch\le c(z,y)\le h
\quad(z\in E,\ y\in S_h^*).
\end{equation}

Take $h_j\rightarrow0$ and John ellipsoids satisfying
\[
y_j+L_{h_j}^*B_1
\subset S_{h_j}^*
\subset y_j+nL_{h_j}^*B_1.
\]
Passing to a subsequence,
\[
y_j\longrightarrow y_\infty,\qquad
L_{h_j}^*\longrightarrow L_\infty^*.
\]
Since $\bigcap_{h>0}\overline{S_h^*}=F$,
\[
y_\infty+L_\infty^*B_1
\subset F
\subset y_\infty+nL_\infty^*\overline B_1.
\]
Therefore
\[
\aff F=y_\infty+\operatorname{range}L_\infty^*,
\]
and
\[
y_\infty+L_\infty^*\overline B_{1/2}\subset\ri F.
\]
This remains valid if $L_\infty^*$ is singular.

Fix a unit vector $e\in V$ and a point $y_*\in\Omega^*$ with $e\cdot y_*\ne0$.
Set
\[
w_j=(L_{h_j}^*)^{-1}(y_*-y_j),
\qquad a_j=\frac{w_j}{|w_j|},
\qquad d_j=L_{h_j}^*a_j.
\]
By \eqref{strip},
\[
|e\cdot y_j|+|(L_{h_j}^*)^Te|\le Ch_j,
\qquad |w_j|\ge\frac1{Ch_j},
\qquad |d_j|\le Ch_j.
\]
For every fixed $R>0$ and all sufficiently large $j$,
\[
B_{1/2}+[0,R]a_j\subset\operatorname{conv}(B_1,w_j),
\]
because $R\le|w_j|/2$.
Convexity of $\Omega^*$ gives
\[
y_j+L_{h_j}^*\bigl(B_{1/2}+[0,R]a_j\bigr)\subset\Omega^*.
\]
The points in the smaller cylinder approach the compact set $y_\infty+L_\infty^*\overline B_{1/4}\subset\ri F$.
By closedness of $\Gamma$ and \eqref{inversecontact},
\[
\begin{aligned}
\delta_j(R)=\sup\{\dist(x,E)\mid {}&(x,y)\in\Gamma,\\
&y\in y_j+L_{h_j}^*(\overline B_{1/4}+[0,R]a_j)\}
\longrightarrow0.
\end{aligned}
\]

Fix $(x,b+td_j)\in\Gamma$, where
\[
b=y_j+L_{h_j}^*z,\qquad |z|\le\frac14,\qquad 0\le t\le R.
\]
Choose $z_0\in E$ such that
\[
D_yc(z_0,b)\cdot d_j=\max_{\zeta\in E}D_yc(\zeta,b)\cdot d_j.
\]
Let
\[
H_x(y)=\frac{c(x,y)-u(x)}{h_j}.
\]
Calibration at $b+td_j$ and the $c$-support with source $z_0$ give
\[
\begin{aligned}
H_x(b)
&\ge\frac{c(z_0,b)}{h_j}
+\frac1{h_j}\int_0^t
\bigl(D_yc(z_0,b+sd_j)-D_yc(x,b+sd_j)\bigr)\cdot d_j\dd s\\
&\ge-C-CR\delta_j(R)-CR^2h_j.
\end{aligned}
\]
For the last inequality, choose $z_1\in E$ with $|x-z_1|\le\delta_j(R)$ and use
\[
\bigl(D_yc(z_0,b)-D_yc(x,b)\bigr)\cdot d_j
\ge-C\delta_j(R)|d_j|,
\]
together with \eqref{strip} and $|d_j|\le Ch_j$.
For fixed $R$, choose $j$ so large that
\[
CR\delta_j(R)+CR^2h_j\le1.
\]
Thus $H_x(b)\ge-C-1$, where $C$ is independent of $R$.
On the base ellipsoid,
\[
H_x(y_j+L_{h_j}^*\zeta)\le\frac{v(y_j+L_{h_j}^*\zeta)}{h_j}\le1
\qquad(\zeta\in B_1).
\]
Since $|(L_{h_j}^*)^{-1}y_j|\le n$, \eqref{propagation} on the segment from $b$ to $0$ yields
\[
H_x(0)=-\frac{u(x)}{h_j}\ge-C_1.
\]
The constant $C_1$ is independent of $R$.
The entire inverse contact image of the smaller cylinder is contained in
\[
\{x\in\overline\Omega\mid u(x)<(C_1+1)h_j\}.
\]
Its intersection with $\Omega$ is contained in $S_{(C_1+1)h_j}$.

Choose an integer $m$ and set $R=m$.
For the pairwise disjoint balls $B_{1/8}(ka_j)$, $k=1,\ldots,m$, we have
\[
\bigcup_{k=1}^mB_{1/8}(ka_j)\subset B_{1/4}+[0,R]a_j.
\]
By mass balance, \eqref{radial} and Proposition~\ref{volume},
\[
\begin{aligned}
C^{-1}m\det L_{h_j}^*
&\le\int_{\bigcup_{k=1}^m\{y_j+L_{h_j}^*B_{1/8}(ka_j)\}}g\dd y\\
&\le\int_{S_{(C_1+1)h_j}}f\dd x\\
&\le C|S_{h_j}|\le C|S_{h_j}^*|\le C\det L_{h_j}^*.
\end{aligned}
\]
Choose $m$ larger than the resulting fixed bound and then let $j$ tend to infinity.
This is a contradiction.
Interchanging the variables proves uniqueness of both closed contact fibers.

Every point of either closed domain has a contact.
Thus the compact graph $\Gamma$ defines inverse continuous mappings.
Invariance of domain, applied to both mappings, gives
\[
T_u(\Omega)=\Omega^*,\qquad T_u(\partial\Omega)=\partial\Omega^*.
\]
Differentiation of the unique maximizing $c$-supports gives the asserted identities for $Du$ and $Dv$.
We henceforth write $T=T_u$ and $T^{-1}=T_v$.

Finally, suppose that sections do not shrink uniformly.
In the original variables, compactness would give
\[
(x_j,y_j)\in\Gamma,\qquad (x_j,y_j)\longrightarrow(x_0,y_0),
\qquad z_j\longrightarrow z_0\ne x_0,
\]
\[
0\le u(z_j)-u(x_j)-c(z_j,y_j)+c(x_j,y_j)\longrightarrow0.
\]
Then $(z_0,y_0)\in\Gamma$, contrary to uniqueness.
The coordinate maps are uniformly bi-Lipschitz, so
\[
\lim_{h\rightarrow 0}\sup_{(x_0,y_0)\in\Gamma}
\bigl(\diam S_h+\diam S_h^*\bigr)=0.
\]
\end{proof}

\section{Normalization of sections}\label{normalization}
In this section, we prove the estimates of sections, which are needed for both Theorem \ref{SR} and \ref{HR}.

\begin{proposition}\label{axes}
For every $\eps>0$, there are constants $c_\eps,C_\eps,h_\eps>0$ such that
\[
 c_\eps h^{1/2+\eps}\le\sigma_{\min}(L_h)
 \le\|L_h\|\le C_\eps h^{1/2-\eps}
 \qquad(0<h<h_\eps).
\]
The same estimates hold for $L_h^*$. The constants are uniform in the contact pair.
\end{proposition}

First we use the normalized coordinates~\eqref{cnorm} and define
\[
 F(h)=\int_{S_h}f\dd x,\qquad
 F^*(h)=\int_{S_h^*}g\dd y,
\]
\[
 h(m)=F^{-1}(m),\qquad h^*(m)=(F^*)^{-1}(m).
\]
For sufficiently small $m>0$, we denote
\[
 I(m)=\int_0^m\bigl[h(s)+h^*(s)\bigr]\dd s,
 \qquad J(m)=m^{-1-2/n}I(m).
\]

\begin{lemma}\label{massmono}
There is a nonnegative function $\omega(m)$, uniform in the contact pair, such that
\begin{equation}\label{monotonicity}
 \omega(m)\longrightarrow0,\qquad
 \frac{d\log J(m)}{d\log m}\ge-\omega(m).
\end{equation}
Moreover,
\begin{equation}\label{heightpower}
 \frac{\log h(m)}{\log m},\quad
 \frac{\log h^*(m)}{\log m}\longrightarrow\frac2n
 \qquad(m\rightarrow 0),
\end{equation}
uniformly in the contact pair.
\end{lemma}

\begin{proof}
Denote $q(t)=c(tx,y)$. It follows by~\eqref{scalar} that
\[
 q''(t)\ge-C|x||y|\,|q'(t)|,\qquad q(0)=0.
\]
If $q(1)\ge0$, by integration of the differential inequality we have
\[
 q(1)\le e^{C|x||y|}q'(1).
\]
Consequently, on the contact set, we obtain
\[
 \begin{aligned}
 x\cdot Du(x)&\ge e^{-C|x||y|}\bigl[u(x)+v(y)\bigr],\\
 y\cdot Dv(y)&\ge e^{-C|x||y|}\bigl[u(x)+v(y)\bigr].
 \end{aligned}
\]
The first expression is positive when $u(x)>0$. Since $Du=0$ almost everywhere on each level set of $u$, we have
\[
 |\{u=t\}|=|\{v=t\}|=0\qquad(t>0).
\]
Thus $F,F^*$ are continuous and strictly increasing. Their inverse functions are also well defined.

It follows from radial integration in the convex source domain that
\[
 \int_{S_h}x\cdot Du\dd x
 \le n\int_{S_h}(h-u)\dd x.
\]
At an endpoint on the source boundary, the omitted boundary term is nonnegative. Freezing $f$ at zero and using the uniform shrinking in Proposition~\ref{contacts}, we obtain
\[
 \int_{\{u<h\}}(u+v)\dd\pi
 \le n(1+o(1))\int_{\{u<h\}}(h-u)\dd\pi.
\]
The error term is uniform in the contact pair. The dual inequality follows by the same argument.

For a measurable set $E$ in the product with $\pi(E)=m$, the layer-cake formula gives
\[
 \int_Ev\dd\pi
 \ge\int_0^\infty\bigl(m-F^*(t)\bigr)_+\dd t
 =\int_0^m h^*(s)\dd s.
\]
Equality holds precisely when $E=\{v<h^*(m)\}$ up to a set of $\pi$-measure zero. The dual case follows by applying this with $E=\{u<h(m)\}$ and interchanging the potentials.

For a common nonnegative error $\omega(m)\to0$, we have
\[
 \begin{aligned}
 \bigl[n(1+\omega(m))+1\bigr]\int_0^m h(s)\dd s
 +\int_0^m h^*(s)\dd s
 &\le n(1+\omega(m))m h(m),\\
 \int_0^m h(s)\dd s
 +\bigl[n(1+\omega(m))+1\bigr]\int_0^m h^*(s)\dd s
 &\le n(1+\omega(m))m h^*(m).
 \end{aligned}
\]
Adding these inequalities implies
\[
 \bigl[n(1+\omega(m))+2\bigr]I(m)
 \le n(1+\omega(m))m I'(m).
\]
Hence we obtain
\[
 \frac{d\log J}{d\log m}
 \ge-\frac{2\omega(m)}{n(1+\omega(m))}
 \ge-\omega(m),
\]
which proves~\eqref{monotonicity}.

It follows from Proposition~\ref{volume} that
\[
 h(m),h^*(m)\ge c m^{2/n},\qquad J(m)\ge c.
\]
For a fixed small $m_0>0$, we have
\[
 J(m)\le J(m_0)
 \exp\left(\int_m^{m_0}\omega(s)\frac{\dd s}{s}\right),
 \qquad
 \int_m^{m_0}\omega(s)\frac{\dd s}{s}=o(|\log m|).
\]

By monotonicity and the preceding bounds,
\[
cm^{2/n}\le h(m),h^*(m)\le\frac{I(2m)}m,
\]
which gives~\eqref{heightpower}. Moreover,
\[
\|Du\|_{L^\infty(\Omega)}+\|Dv\|_{L^\infty(\Omega^*)}\le C,
\qquad
|\Omega\cap B_r|,|\Omega^*\cap B_r|\ge cr^n
\quad(0<r<r_0).
\]
Hence $F(t),F^*(t)\ge ct^n$ for $0<t<t_0$. The constants and the initial mass $m_0$ are uniform in the contact pair.
\end{proof}

Set $a_n=(n+2)/(2n)$.

\begin{lemma}\label{homogeneous}
Suppose $m_j\rightarrow0$ and
\[
\frac{J(m_jt)}{J(m_j)}\longrightarrow1
\quad\text{locally uniformly for }t>0.
\]
Let
\[
\ell_j=\frac{I(m_j)}{m_j},
\qquad L_j=L_{h(m_j)},
\qquad L_j^*=L_{h^*(m_j)}^*,
\]
\[
\tilde u_j(X)=\ell_j^{-1}u(L_jX),
\qquad
\tilde v_j(Y)=\ell_j^{-1}v(L_j^*Y).
\]
After a subsequence,
\[
\tilde u_j\longrightarrow\tilde u\quad\text{in }C_{\rm loc}(\tilde\Omega),
\qquad
\tilde v_j\longrightarrow\tilde v\quad\text{in }C_{\rm loc}(\tilde\Omega^*),
\]
where $\tilde\Omega,\tilde\Omega^*$ are open convex cones and $\tilde u,\tilde v$ are finite convex functions satisfying
\[
\tilde u(rX)=r^2\tilde u(X),
\qquad
\tilde v(rY)=r^2\tilde v(Y)
\qquad(r>0).
\]
In particular,
\[
\begin{aligned}
K_t&=\overline{\{\tilde u<t\}}=\sqrt{t/a_n}\,K_{a_n},\\
K_t^*&=\overline{\{\tilde v<t\}}=\sqrt{t/a_n}\,K_{a_n}^*
\qquad(t>0).
\end{aligned}
\]
The base contact pair may vary with $j$.
\end{lemma}

\begin{proof}
By the definition of $\ell_j$,
\[
\int_0^t\frac{h(m_js)+h^*(m_js)}{\ell_j}\,\dd s
=\frac{I(m_jt)}{I(m_j)}
\longrightarrow t^{1+2/n}.
\]
The two summands are increasing and locally bounded, with
\[
\sup_j\int_0^\delta\frac{h(m_js)+h^*(m_js)}{\ell_j}\,\dd s
\le C\delta
\qquad(0<\delta<1).
\]
Pass to subsequential limits in the two inequalities preceding~\eqref{monotonicity}. Their sum is an equality. If $d$ is the difference of the two limits, subtraction gives
\[
\int_0^t d(s)\,\dd s=td(t),
\qquad
\left|\int_0^t d(s)\,\dd s\right|\le t^{1+2/n}.
\]
Thus $d'=0$ in distributions and $|d(t)|\le t^{2/n}$, so $d=0$. Consequently,
\begin{equation}\label{inverseheights}
\frac{h(m_jt)}{\ell_j},
\qquad
\frac{h^*(m_jt)}{\ell_j}
\longrightarrow a_nt^{2/n}
\end{equation}
locally uniformly for $t>0$.

Set
\[
\tilde\pi_j=m_j^{-1}\bigl(L_j^{-1},(L_j^*)^{-1}\bigr)_\#\pi.
\]
By~\eqref{mass}, John's lemma and Proposition~\ref{contacts},
\[
c\le\frac{|\det L_j|}{m_j},\frac{|\det L_j^*|}{m_j}\le C,
\qquad
\|L_j\|+\|L_j^*\|\longrightarrow0.
\]
After a subsequence,
\[
\frac{|\det L_j|}{m_j}f(0)\longrightarrow f_0>0,
\qquad
\frac{|\det L_j^*|}{m_j}g(0)\longrightarrow g_0>0.
\]
For every fixed $R>0$, continuity gives
\[
\begin{aligned}
&\sup_{X\in L_j^{-1}\Omega\cap B_R}
\left|\frac{|\det L_j|}{m_j}f(L_jX)-f_0\right|\\
&\quad+\sup_{Y\in(L_j^*)^{-1}\Omega^*\cap B_R}
\left|\frac{|\det L_j^*|}{m_j}g(L_j^*Y)-g_0\right|
\longrightarrow0.
\end{aligned}
\]
By~\eqref{inverseheights},
\[
\frac{F(t\ell_j)}{m_j},
\qquad
\frac{F^*(t\ell_j)}{m_j}
\longrightarrow(t/a_n)^{n/2}.
\]
For each $t>0$, radial comparison and the volume bounds give
\[
L_j^{-1}S_{t\ell_j}\subset B_{C_t},
\qquad
c_t\le|L_j^{-1}S_{t\ell_j}|\le C_t,
\]
and the same estimates hold on the target side. A diagonal subsequence gives Hausdorff limits
\[
\overline{L_j^{-1}S_{t\ell_j}}\longrightarrow K_t,
\qquad
\overline{(L_j^*)^{-1}S_{t\ell_j}^*}\longrightarrow K_t^*
\qquad(t\in\mathbb Q_+),
\]
with
\[
f_0|K_t|=g_0|K_t^*|=(t/a_n)^{n/2}.
\]
Nesting and continuity of the volumes extend these limits to every $t>0$.

Let
\[
\tilde\Omega=\bigcup_{t>0}\operatorname{int}K_t,
\qquad
\tilde\Omega^*=\bigcup_{t>0}\operatorname{int}K_t^*.
\]
These sets are open and convex. By~\eqref{propagation} and~\eqref{rescaledscalar}, for every $K\Subset\tilde\Omega$,
\[
0\le\tilde u_j\le C_K,
\qquad
|D\tilde u_j|\le C_K,
\qquad
D^2\tilde u_j+C_K\|L_j\|\operatorname{Id}\ge0
\quad\text{on }K.
\]
The last inequality holds in distributions by semiconvexity. Hence, after a subsequence,
\[
\tilde u_j\longrightarrow\tilde u
\quad\text{locally uniformly},
\qquad
D\tilde u_j\longrightarrow D\tilde u
\quad\text{in }L^1_{\rm loc}\text{ and a.e.}
\]
The limit is finite and convex. The same conclusions hold for $\tilde v_j$. Radial comparison gives
\[
\tilde u(0)=\tilde v(0)=0,
\qquad
X\cdot D\tilde u\ge\tilde u,
\qquad
Y\cdot D\tilde v\ge\tilde v.
\]
For $t>0$,
\[
|\{\tilde u=t\}|=|\{\tilde v=t\}|=0,
\qquad
\{\tilde u<t\}\subset\operatorname{int}K_t\subset\{\tilde u\le t\}.
\]
Convexity therefore yields
\[
K_t=\overline{\{\tilde u<t\}},
\qquad
K_t^*=\overline{\{\tilde v<t\}}.
\]

The radial integral inequality gives
\[
\int_{\{u<h\}}v\,\dd\pi\le ChF(h).
\]
Thus, for each $H>0$,
\[
\tilde\pi_j\{\tilde u_j<H,\ \tilde v_j>R\}
+\tilde\pi_j\{\tilde v_j<H,\ \tilde u_j>R\}
\le\frac{C_H}{R}.
\]
Together with the boundedness of sections, this gives a locally weakly convergent subsequence $\tilde\pi_j\rightharpoonup\tilde\pi$ on $\tilde\Omega\times\tilde\Omega^*$ retaining both marginals,
\begin{equation}\label{limitmarginals}
(\tilde\pi)_X=f_0\mathbf1_{\tilde\Omega}\,\dd X,
\qquad
(\tilde\pi)_Y=g_0\mathbf1_{\tilde\Omega^*}\,\dd Y.
\end{equation}
In particular,
\[
\tilde\pi\{\tilde u<t\}
=\tilde\pi\{\tilde v<t\}
=(t/a_n)^{n/2}.
\]

Let $0\le\varphi\in C_c(\tilde\Omega\times\tilde\Omega^*)$ and choose $K\Subset\tilde\Omega$ containing the source projection of its support. The marginal bounds give
\[
\left|\int\varphi X\cdot(D\tilde u_j-D\tilde u)\,\dd\tilde\pi_j\right|
\le C_\varphi\|D\tilde u_j-D\tilde u\|_{L^1(K)}
\longrightarrow0.
\]
Approximation by continuous functions and local weak convergence then yield
\[
\int\varphi X\cdot D\tilde u_j\,\dd\tilde\pi_j
\longrightarrow
\int\varphi X\cdot D\tilde u\,\dd\tilde\pi.
\]
Since $|L_jX||L_j^*Y|\to0$ on $\operatorname{spt}\varphi$, the radial inequalities pass to the limit,
\[
\begin{aligned}
X\cdot D\tilde u(X)&\ge\tilde u(X)+\tilde v(Y),\\
Y\cdot D\tilde v(Y)&\ge\tilde u(X)+\tilde v(Y)
\qquad\tilde\pi\text{-a.e.}
\end{aligned}
\]

For each direction meeting $\tilde\Omega$, let
\[
r_t(\theta)=\sup\{r>0\mid r\theta\in\tilde\Omega,\ \tilde u(r\theta)<t\},
\qquad
b_t(\theta)=\lim_{r\uparrow r_t(\theta)}\tilde u(r\theta)\le t.
\]
Integration along rays gives
\[
\begin{aligned}
n\int_{\{\tilde u<t\}}(t-\tilde u)\,\dd X
-\int_{\{\tilde u<t\}}X\cdot D\tilde u\,\dd X
&=\int r_t(\theta)^n[t-b_t(\theta)]\,\dd\theta\\
&\ge0.
\end{aligned}
\]
The angular integral is over the directions meeting $\tilde\Omega$. By the common distribution function,
\[
\int_{\{\tilde u<t\}}\tilde u\,\dd\tilde\pi
=\int_{\{\tilde v<t\}}\tilde v\,\dd\tilde\pi
=\frac{n}{n+2}t(t/a_n)^{n/2}.
\]
Rearrangement and radial integration give
\[
\begin{aligned}
2\int_{\{\tilde u<t\}}\tilde u\,\dd\tilde\pi
&\le\int_{\{\tilde u<t\}}(\tilde u+\tilde v)\,\dd\tilde\pi\\
&\le f_0\int_{\{\tilde u<t\}}X\cdot D\tilde u\,\dd X\\
&\le n\int_{\{\tilde u<t\}}(t-\tilde u)\,\dd\tilde\pi
=2\int_{\{\tilde u<t\}}\tilde u\,\dd\tilde\pi.
\end{aligned}
\]
All inequalities are equalities. Equality in rearrangement at rational heights gives
\[
\tilde u(X)=\tilde v(Y)\quad\tilde\pi\text{-a.e.}
\]
Using~\eqref{limitmarginals} and the dual argument,
\[
X\cdot D\tilde u=2\tilde u,
\qquad
Y\cdot D\tilde v=2\tilde v.
\]
Consequently,
\[
\begin{aligned}
\tilde u(rX)&=r^2\tilde u(X)
&&\text{if }X,rX\in\tilde\Omega,\\
\tilde v(rY)&=r^2\tilde v(Y)
&&\text{if }Y,rY\in\tilde\Omega^*.
\end{aligned}
\]

Write $r_{\max}(\theta)=\sup\{r>0\mid r\theta\in\tilde\Omega\}$. If $r_{\max}(\theta)<\infty$, homogeneity gives a finite limit of $\tilde u(r\theta)$ as $r\uparrow r_{\max}(\theta)$. For some $t\in\mathbb Q_+$,
\[
r_t(\theta)=r_{\max}(\theta),
\qquad
r_t(\theta)^n[t-b_t(\theta)]>0.
\]
Equality in radial integration at all rational heights therefore gives $r_{\max}(\theta)=\infty$ for almost every direction meeting $\tilde\Omega$.

Suppose $R\theta\in\partial\tilde\Omega$ and $r_0\theta\in\tilde\Omega$ for some $0<r_0<R$. A supporting hyperplane gives
\[
\xi\cdot X<s\quad(X\in\tilde\Omega),
\qquad
R\xi\cdot\theta=s>0.
\]
For all $\theta'$ in a neighborhood of $\theta$,
\[
r_0\theta'\in\tilde\Omega,
\qquad
\xi\cdot\theta'>0,
\qquad
r_{\max}(\theta')\le\frac{s}{\xi\cdot\theta'}<\infty,
\]
a contradiction. Hence
\[
r\tilde\Omega=\tilde\Omega,
\qquad
r\tilde\Omega^*=\tilde\Omega^*
\qquad(r>0).
\]
The section identities follow. All estimates are uniform in the base contact pair.
\end{proof}

\begin{proof}[Proof of Proposition~\ref{axes}]
For $0<a<b<m_0$, set
\[
\mathcal D(a,b)=\log\frac{J(b)}{J(a)}
+\int_a^b\omega(s)\frac{\dd s}{s}\ge0.
\]
Then
\[
\mathcal D(a,b')=\mathcal D(a,b)+\mathcal D(b,b')
\qquad(0<a<b<b'<m_0).
\]
Fix $0<\theta<1$ and put $m_k=\theta^km_0$. Since $J\ge C^{-1}$ and $\omega(m)\to0$,
\[
\begin{aligned}
\sum_{k=0}^{N-1}\mathcal D(m_{k+1},m_k)
&=\log\frac{J(m_0)}{J(m_N)}
+\int_{m_N}^{m_0}\omega(s)\frac{\dd s}{s}\\
&\le C+\int_{m_N}^{m_0}\omega(s)\frac{\dd s}{s}
=o(N).
\end{aligned}
\]

Write $L(m)=L_{h(m)}$. We claim that for every $\eta>0$ there are $N_0\in\mathbb N$ and $\tau,m_*>0$, with $m_*\theta^{-N_0}<m_0$, such that
\[
0<m<m_*,
\qquad
\mathcal D(m\theta^{N_0},m\theta^{-N_0})\le\tau
\]
imply
\begin{equation}\label{johnstep}
\begin{aligned}
e^{-\eta}\theta^{2/n}L(m)L(m)^T
&\le L(\theta m)L(\theta m)^T\\
&\le e^\eta\theta^{2/n}L(m)L(m)^T.
\end{aligned}
\end{equation}
Otherwise, there are $m_j\rightarrow0$ such that
\[
\mathcal D(m_j\theta^j,m_j\theta^{-j})\le j^{-1}
\]
and~\eqref{johnstep} fails at $m_j$. Additivity and $\omega(m)\to0$ give
\[
\frac{J(m_jt)}{J(m_j)}\longrightarrow1
\quad\text{locally uniformly for }t>0.
\]
Lemma~\ref{homogeneous} and continuity of John ellipsoids imply
\[
L(m_j)^{-1}L(\theta m_j)L(\theta m_j)^TL(m_j)^{-T}
\longrightarrow\theta^{2/n}\operatorname{Id},
\]
contrary to the failure of~\eqref{johnstep}.

By additivity,
\[
\begin{aligned}
&\#\{N_0\le k<N-N_0\mid
\mathcal D(m_{k+N_0},m_{k-N_0})>\tau\}\\
&\qquad\le\frac{2N_0}{\tau}
\sum_{k=0}^{N-1}\mathcal D(m_{k+1},m_k)
=o(N).
\end{aligned}
\]
Thus~\eqref{johnstep} fails at only $o(N)$ of the first $N$ scales.

Set $L_k=L(m_k)$. Nesting and the fixed mass ratio give
\[
\|L_k^{-1}L_{k+1}\|\le n,
\qquad
|\det(L_k^{-1}L_{k+1})|\ge c\theta.
\]
Hence, for some $c_0>0$,
\[
c_0^2L_kL_k^T\le L_{k+1}L_{k+1}^T\le n^2L_kL_k^T.
\]
Combining this with~\eqref{johnstep},
\[
\begin{aligned}
e^{-\eta N-o(N)}\theta^{2N/n}L_0L_0^T
&\le L_NL_N^T\\
&\le e^{\eta N+o(N)}\theta^{2N/n}L_0L_0^T.
\end{aligned}
\]
Since $C^{-1}\operatorname{Id}\le L_0L_0^T\le C\operatorname{Id}$, interpolation between consecutive mass scales and then $\eta\rightarrow0$ give
\[
\frac{\log\sigma_{\min}(L(m))}{\log m},
\qquad
\frac{\log\|L(m)\|}{\log m}
\longrightarrow\frac1n.
\]
The assertion follows from~\eqref{heightpower}. The proof for $L_h^*$ is identical.
\end{proof}

\begin{corollary}\label{growth}
Fix $0<\eps<1/4$ and set
\[
b_-=\frac12-\eps,
\qquad
b_+=\frac12+\eps,
\qquad
a=\frac{b_-}{b_+},
\qquad
q=\frac1{b_-}.
\]
In the coordinates~\eqref{cnorm},
\[
\begin{aligned}
c_\eps|x|^q&\le u(x)\le C_\eps|x|^{1+a},\\
c_\eps|y|^q&\le v(y)\le C_\eps|y|^{1+a},
\end{aligned}
\]
\[
\Omega\cap B_{c_\eps h^{b_+}}
\subset S_h
\subset\Omega\cap B_{C_\eps h^{b_-}},
\]
with the dual inclusions for $S_h^*$. Moreover,
\[
|T(x)-T(z)|+|Du(x)-Du(z)|
\le C_\eps|x-z|^a
\qquad(x,z\in\overline\Omega),
\]
\[
|T^{-1}(y)-T^{-1}(w)|+|Dv(y)-Dv(w)|
\le C_\eps|y-w|^a
\qquad(y,w\in\overline{\Omega^*}).
\]
The constants are uniform in the contact pair.
\end{corollary}

\begin{proof}
By Proposition~\ref{axes} and $0\in\overline{S_h}\cap\overline{S_h^*}$,
\[
S_h,S_h^*\subset B_{C_\eps h^{b_-}}.
\]
Hence $u(x)\ge c_\eps|x|^q$ and $v(y)\ge c_\eps|y|^q$. For $y=T(x)$,
\[
c_\eps|y|^q\le v(y)\le u(x)+v(y)=c(x,y)\le C|x||y|.
\]
Since $1/(q-1)=a$,
\[
|T(x)|\le C_\eps|x|^a,
\qquad
|Du(x)|=|c_x(x,T(x))|\le C_\eps|x|^a.
\]
Integration along $[0,x]$ gives
\[
u(x)\le C_\eps|x|^{1+a},
\qquad
\frac1{1+a}=b_+,
\]
and therefore the section inclusions. Recenter at each contact pair and use the uniform bounds for the coordinate changes to obtain the estimates for $T$ and $Du$. Interchanging source and target gives the dual estimates.
\end{proof}

Put
\[
P_h=hL_h^{-T},
\qquad
\tilde c_h(X,Y)=h^{-1}c(L_hX,P_hY).
\]
Proposition~\ref{axes} gives
\[
\begin{aligned}
\|L_h\|+\|P_h\|&\le C_\eps h^{1/2-\eps},\\
\|L_h^{-1}\|+\|P_h^{-1}\|&\le C_\eps h^{-1/2-\eps}.
\end{aligned}
\]
The coordinate identities~\eqref{cnorm} and the $C^4$ bound for the original cost imply
\[
|D_x^iD_y^j(c-x\cdot y)|
\le C|x|^{2-i}|y|^{2-j}
\qquad(0\le i,j\le2,\ i+j\le3),
\]
\[
|D_x^3c|\le C|y|,
\qquad
|D_y^3c|\le C|x|.
\]
Thus, for every fixed $R>0$ and $0<\eps<1/4$,
\begin{equation}\label{costlimit}
\begin{aligned}
\|\tilde c_h-X\cdot Y\|_{C^{2,1}(B_R\times B_R)}
&\le C_Rh^{-1}(\|L_h\|+\|P_h\|)^4\\
&\le C_{\eps,R}h^{1-4\eps}
\longrightarrow0.
\end{aligned}
\end{equation}

\section{Obliqueness}\label{obliqueness}

We apply Corollary~\ref{growth} to prove the obliqueness estimate used in Theorem~\ref{HR}.

\begin{proposition}\label{oblique}
Assume $\partial\Omega,\partial\Omega^*\in C^{1,\beta}$ for some $0<\beta\le1$.

There exists $\theta_0>0$ such that
\[
\nu_\Omega(x)^T c_{xy}(x,T(x))^{-T}\nu_{\Omega^*}(T(x))\ge\theta_0
\qquad(x\in\partial\Omega),
\]
where $\nu_\Omega$ and $\nu_{\Omega^*}$ are the inward unit normals.
\end{proposition}

\begin{proof}
For convenience, we assume that $(0,0)$ is a boundary contact pair . By the coordinates in \eqref{cnorm}, we also assume that both domains are convex. Denote their inward normals at the origin by $\nu$ and $\nu^*$.

Choose $0<\varepsilon<\min\{\beta,1\}/100$. Then it follows from Corollary~\ref{growth} that
\[
    u(x)\ge c_\varepsilon |x|^q,\qquad |T(x)|\le C_\varepsilon |x|^a.
\]
Suppose that $\nu\cdot\nu^*\le0$.
Near the origin, we assume
\[
\Omega=\{z+s\nu\mid s>\phi(z),\ z\cdot\nu=0\},\qquad 0\le\phi(z)\le K|z|^{1+\beta}.
\]
For $K_1>K$ and sufficiently small $t>0$, we define
\[
x_t=-t\nu^*+K_1t^{1+\beta}\nu,\qquad y_t=T(x_t).
\]
Then it follows that
\[
x_t\cdot\nu\ge K_1t^{1+\beta},\qquad |x_t-(x_t\cdot\nu)\nu|\le t,\qquad x_t\in\Omega,\qquad t\le |x_t|\le 2t.
\]
By convexity of $\Omega^*$ we have $\nu^*\cdot y_t\ge0$.
Since $u,v\ge0$ and $u(x_t)+v(y_t)=c(x_t,y_t)$, we obtain that
\[
    c_\varepsilon t^q
    \le u(x_t)\le c(x_t,y_t) \le C\bigl(t^{1+\beta}|y_t|+t^2|y_t|^2\bigr) \le C_\varepsilon\bigl(t^{1+\beta+a}+t^{2+2a}\bigr).
\]
But we note that
\[
1+\beta+a>q,\qquad 2+2a>q.
\]
Letting $t\to0$ proves $\nu\cdot\nu^*>0$.

At every boundary contact pair, we can use the same arguments to obtain the asserted inequality. Finally, a uniform positive constant $\theta_0$ follows from Proposition~\ref{contacts} and compactness arguments.
\end{proof}

\section{Global \texorpdfstring{$W^{2,p}$}{W2,p} estimates}\label{sobolev}

In this section, we first prove an estimate on interior sections. All constants are uniform as the center of section approaches the boundary.

\begin{lemma}\label{interior}
Let $1<p<\infty$ and $0<\eps<1/10$. In the coordinates~\eqref{cnorm}, suppose $\overline{S_h}\subset\Omega$. For all sufficiently small $h>0$,
\[
\int_{S_{h/2}}|D^2u|^p\dd x\le C_{p,\eps}h^{n/2-2\eps p}.
\]
\end{lemma}

\begin{proof}
Define
\[
D=L_h^{-1}S_h,\qquad P_h=hL_h^{-T},\qquad
\tilde u(X)=h^{-1}u(L_hX),\qquad
\tilde c(X,Y)=h^{-1}c(L_hX,P_hY).
\]
Then there exists a positive constant $R$ depending only on $n$ such that
\[
B_1(X_h)\subset D\subset B_R,\qquad
\tilde u(0)=0,\qquad \tilde u=1\quad\text{on }\partial D.
\]
For $|q|<1/(4R)$, we obtain by \eqref{costlimit} that
\[
\tilde u-\tilde c(\cdot,q)\ge\frac12\quad\text{on }\partial D,
\qquad
\bigl(\tilde u-\tilde c(\cdot,q)\bigr)(0)=0.
\]
Hence the minimum of $\tilde u-\tilde c(\cdot,q)$ on $\overline{D}$ is reached at a interior point $X$, which yields that
\[
D\tilde u(X)=D_X\tilde c(X,q).
\]
Moreover, we have
\[
P_hq=T(L_hX),\qquad
B_{1/(4R)}\subset P_h^{-1}T(S_h).
\]
It follows from mass balance and Proposition~\ref{volume} that
\begin{equation}\label{interiorvolume}
|\det P_h|\le C|\det L_h|,\qquad
c\le\frac{(\det L_h)^2}{h^n}\le C.
\end{equation}
By Corollary~\ref{growth} and $c_{xx}(x,y)=O(|y|^2)$ we obtain, on the actual transport graph,
\[
\mathcal A_h(X):=h^{-1}L_h^Tc_{xx}(L_hX,T(L_hX))L_h,
\]
\[
\|\mathcal A_h\|_{L^\infty(D)}\le C h^{2b_-^2/b_+-2\eps}=\epsilon_h\longrightarrow0.
\]
In particular, $\tilde u+\epsilon_h|X|^2/2$ is convex.

The normalized equation is
\begin{equation}\label{normalizedequation}
\det(D^2\tilde u-\mathcal A_h)=\kappa_h(1+r_h),\qquad
\kappa_h=\frac{(\det L_h)^2}{h^n}\frac{f(0)}{g(0)},
\end{equation}
where
\[
c\le\kappa_h\le C,\qquad \|r_h\|_{L^\infty(D)}=\delta_h\longrightarrow0.
\]
The factor $r_h$ includes the oscillations of $f$, $g$ and $\det c_{xy}$.
By Jiang~\cite[Theorem 1.1(i)--(ii)]{Jiang2023}, we deduce that the equation holds in the admissible viscosity sense.

Let $z_D$ be the convex solution of
\[
\det D^2z_D=1\quad\text{in }D,\qquad z_D=0\quad\text{on }\partial D.
\]
Existence and uniqueness follow from~\cite[Theorem 2.1]{TW2008}. Define
\[
w=1+\kappa_h^{1/n}z_D.
\]
For $h>0$ small, we consider two test functions
\[
\begin{aligned}
w_-&=1+[\kappa_h(1+2\delta_h)]^{1/n}z_D
       +2\epsilon_h(|X|^2-R^2),\\
w_+&=1+[\kappa_h(1-2\delta_h)]^{1/n}z_D
       -2\epsilon_h(|X|^2-R^2).
\end{aligned}
\]
At an interior contact, the lower test has augmented determinant larger than the right hand side of~\eqref{normalizedequation}. An admissible upper test has augmented determinant smaller than that right hand side. The boundary inequalities and viscosity comparison yield
\begin{equation}\label{interiorcomparison}
w_-\le\tilde u\le w_+,\qquad
\|\tilde u-w\|_{L^\infty(D)}\le C(\epsilon_h+\delta_h)\longrightarrow0.
\end{equation}

The Aleksandrov boundary estimate~\cite[Lemma 3.3]{TW2008} gives
\[
|z_D(X)|\le C\dist(X,\partial D)^{1/n}.
\]
The Pogorelov estimate and interior regularity~\cite[Lemma 3.6 and Theorem 3.1]{TW2008}, applied first on $\{w<7/8\}$, imply
\begin{equation}\label{modelbounds}
\dist(\{w\le3/4\},\partial D)\ge c,\qquad
cI\le D^2w\le CI,\qquad \|w\|_{C^3(\{w\le3/4\})}\le C.
\end{equation}
It follows from semiconvexity and~\eqref{interiorcomparison} that
\[
\|D\tilde u-Dw\|_{L^\infty(\{w\le2/3\})}\longrightarrow0.
\]
For $Y=P_h^{-1}T(L_hX)$, we have
\[
D_X\tilde c(X,Y)-Y
=\int_0^1D^2_{XX}\tilde c(tX,Y)X\dd t
=O(\epsilon_h).
\]
Consequently, for $\tilde T:=P_h^{-1}T(L_h\cdot)$, we obtain
\[
\|\tilde T-Dw\|_{L^\infty(\{w\le2/3\})}\longrightarrow0.
\]

Fix $X_0\in\{\tilde u\le1/2\}$. Uniform convexity in~\eqref{modelbounds}, interior minimization and injectivity give fixed radii $r,r',\rho>0$ such that
\[
B_\rho(\tilde T(X_0))\subset\tilde T(B_{2r}(X_0)),\qquad
\tilde T(B_{r'}(X_0))\subset B_\rho(\tilde T(X_0)).
\]
The full preimage of this target ball is contained in $B_{2r}(X_0)$. We restrict the transport plan to this preimage and the target ball.

Denote $\lambda_h=\kappa_h^{1/n}$ and $Y_0=\tilde T(X_0)$. Under
\[
X=X_0+\widehat X,\qquad Y=Y_0+\lambda_h\widehat Y,
\]
consider the cost
\[
\widehat c(\widehat X,\widehat Y)
=\lambda_h^{-1}\bigl[\tilde c(X,Y)-\tilde c(X,Y_0)-\tilde c(X_0,Y)+\tilde c(X_0,Y_0)\bigr].
\]
Dividing both marginal measures by $|\det L_h|f(0)$, we obtain
\[
\widehat f(\widehat X)=\frac{f(L_hX)}{f(0)},\qquad
\widehat g(\widehat Y)=\frac{g(P_hY)}{g(0)}.
\]
After this bounded change of variables, the two supports $C_1,C_2$ satisfy
\[
B_{1/K}\subset C_1,C_2\subset B_K,\qquad C_2\text{ is convex},
\]
\[
\|\widehat f-1\|_\infty+\|\widehat g-1\|_\infty=o(1),\qquad
\|\widehat c-\widehat X\cdot\widehat Y\|_{C^{1,1}}=o(1),\qquad
\|\widehat c\|_{C^{2,1}}\le C.
\]
By taking its $c$-transform over $C_2$, we extends the potential to $B_K$.

Now we apply the results of Chen and Figalli in~\cite[Proposition 2.4]{CF2017} to the minimizing cost $-\widehat c$ and obtain a local $W^{2,p}$ estimate. Furthermore, by a fixed finite covering of $\{\tilde u\le1/2\}$, we obtain
\[
\int_{\{\tilde u\le1/2\}}|D^2\tilde u|^p\dd X\le C_p.
\]
By Proposition~\ref{axes} and~\eqref{interiorvolume},
\[
h\|L_h^{-1}\|^2\le C_\eps h^{-2\eps},\qquad
|\det L_h|\le Ch^{n/2}.
\]
Scaling back proves the lemma. The changes of $c$-exponential coordinates contribute only bounded lower order terms.
\end{proof}

\begin{corollary}\label{interiorholder}
If $f,g\in C^{\alpha}$ with $0<\alpha<1$, then $u\in C^{2,\alpha}_{\mathrm{loc}}(\Omega)$ and $v\in C^{2,\alpha}_{\mathrm{loc}}(\Omega^*)$.
\end{corollary}

\begin{proof}
Use the same notations in~Lemma \ref{interior}. For $X_0\in\{\tilde u\le1/2\}$, we denote $H=D^2w(X_0)$ and $Y_0=\tilde T(X_0)$. Consider the following change
\[
X=X_0+rH^{-1/2}\xi,\qquad Y=Y_0+rH^{1/2}\eta.
\]
Then we subtract the $c$-support at $(X_0,Y_0)$ and divide the potential and cost by $r^2$. Since $\det H=\kappa_h$, a common mass normalization again gives source and target densities $f/f(0)$ and $g/g(0)$.

First fix $r>0$ small and then take $h>0$ small. It follows that the new potential is arbitrarily close to $|\xi|^2/2$ on $B_2$, the cost is close to $\xi\cdot\eta$, and the new map is close to the identity. As in \cite[Section 3, (3.9)--(3.10)]{DPF2015}, we consider the restricted full plan on
\[
C_1=\overline B_1,\qquad
C_2=r^{-1}H^{-1/2}\bigl(\tilde T(X_0+rH^{-1/2}C_1)-Y_0\bigr),
\]
\[
B_{1/3}\subset C_1,C_2\subset B_3.
\]
The last inclusion follows from interior minimization and injectivity.  Applying \cite[Theorem 5.3]{DPF2015} to the negative of the normalized cost, we prove the $C^{2,\alpha}$ regularity.
\end{proof}

\begin{proof}[Proof of Theorem~\ref{SR}]
A fixed $c$-exponential map sends $\Omega$ to a bounded convex domain and is bi-Lipschitz on the closure. Hence
\[
\bigl|\{x\in\Omega\mid\dist(x,\partial\Omega)<d\}\bigr|\le Cd.
\]
Let $h=c d^{1/b_-}$ and $r=c'h^{b_+}$.
For any $x$ in the set $\{d/2<\dist(x,\partial\Omega)<d\}$, it follows from
Corollary~\ref{growth} that
\[
\overline{S_h(x)}\subset\Omega,\qquad B_r(x)\subset S_{h/2}(x).
\]
A disjoint family of balls $B_{r/5}$ whose corresponding balls $B_r$ cover $\{d/2<\dist(x,\partial\Omega)<d\}$ has at most $Cdr^{-n}$ members. By Lemma~\ref{interior}, we obtain
\[
\begin{aligned}
\int_{\{d/2<\dist(x,\partial\Omega)<d\}}|D^2u|^p\dd x
&\le Cdr^{-n}h^{n/2-2\eps p}\\
&\le C d^{\,1-\eps(n+2p)/(1/2-\eps)}.
\end{aligned}
\]
Choose
\[
0<\eps<\min\left\{\frac1{10},\frac1{2(n+2p+1)}\right\}.
\]
Summing over $d=2^{-k}d_0$ and covering the compact interior, we have
\[
\int_\Omega|D^2u|^p\dd x\le C_p.
\]
The local $W^{2,p}$ regularity identifies the distributional Hessian. Since $u$ and $Du$ are bounded, $u\in W^{2,p}(\Omega)$. The same argument applies to $v$.

The homeomorphism assertion follows from Proposition~\ref{contacts}. Since
\[
\frac{1/2-\eps}{1/2+\eps}\longrightarrow1\qquad\text{as }\eps\longrightarrow0,
\]
Corollary~\ref{growth} gives every H\"older exponent $\eta<1$ for both maps.
\end{proof}

\section{Global H\"older estimates}\label{holderregularity}
Throughout this section, we fix
\[
0<\varepsilon<\frac1{100}\min\{\alpha,\beta\},\qquad b_-=\frac12-\varepsilon,\qquad b_+=\frac12+\varepsilon,\qquad a=\frac{b_-}{b_+},\qquad q=\frac1{b_-}.
\]

\begin{lemma}\label{boundarycomparison}
In the coordinates \eqref{cnorm} at a boundary contact pair, there is a fixed function $\varphi\in C^{2,1}$ such that $\hat u=u+\varphi$ is convex and nondecreasing in the common normal direction.
There are bounded convex domains $D_h$ and convex solutions $w_h$ with
\[
\det D^2w_h=\frac{f(0)}{g(0)}\quad\hbox{in }D_h,\qquad w_h=h\quad\hbox{on }\partial D_h,
\]
\[
B_{c h^{b_+}}\subset D_h\subset B_{C h^{b_-}},\qquad c h^{n/2}\le |D_h|\le C h^{n/2},\qquad D_{h/4}\subset D_h.
\]
For
\[
\delta=\min\bigl\{a\alpha b_-,\ b_-\bigl[1+a(1+\beta)\bigr]-1\bigr\}>2\varepsilon,
\]
we have
\begin{equation}\label{modelerror}
\|\hat u-w_h\|_{L^\infty(D_h\cap\Omega)}+\|w_h-w_{h/4}\|_{L^\infty(D_{h/4})}\le C h^{1+\delta}.
\end{equation}
\end{lemma}

\begin{proof}
By Proposition~\ref{oblique}, a uniformly bounded slide transform makes both inward normals $e_n$.
The normalization \eqref{cnorm} is preserved.
Near the origin, we have
\[
\Omega\subset\{x_n\ge0\},\qquad \Omega^*\subset\{y_n\ge0\},
\]
and assume that $\partial\Omega,\partial\Omega^*$ are represented by
\[
x_n=\phi(x'),\qquad y_n=\phi^*(y'),\qquad 0\le\phi(x'),\phi^*(x')\le C|x'|^{1+\beta}.
\]
It follows from Corollary~\ref{growth} that
\begin{equation}\label{initialgrowth}
c|x|^q\le u(x)\le C|x|^{1+a},\qquad |T(x)|+|Du(x)|\le C|x|^a,
\end{equation}
\[
|T(x)-T(z)|+|Du(x)-Du(z)|\le C|x-z|^a.
\]
Note that the normalized cost satisfies
\[
|c_x(x,y)-y|\le C|x||y|^2,\qquad |c_{xx}(x,y)|\le C|y|^2.
\]
Combined with Corollary~\ref{interiorholder}, we have
\[
D^2u\ge-C|x|^{2a}I,\qquad u_n\ge-C|x|^{1+2a}.
\]

Choose large constants $C_0,K$ and define
\begin{equation}\label{correction}
\varphi(x)=K\bigl(C_0|x|^{2+2a}+x_n|x|^{1+2a}\bigr),\qquad \hat u=u+\varphi.
\end{equation}
For $x_n\ge0$, we have
\[
    K^{-1}\varphi_n
    =C_0(2+2a)x_n|x|^{2a}+|x|^{1+2a}+(1+2a)x_n^2|x|^{2a-1} \ge |x|^{1+2a}.
\]
Taking $C_0$ first and then $K$ sufficiently large, we have
\[
D^2\varphi\ge K|x|^{2a}I,\qquad \varphi_n\ge K|x|^{1+2a},\qquad |D^j\varphi(x)|\le C|x|^{2+2a-j}\quad(0\le j\le3).
\]
Thus
\[
D^2\hat u\ge0,\qquad \hat u_n\ge0,\qquad \varphi\in C^{2,1}.
\]

Extend $\hat u$ by the supremum of its supporting planes based in a fixed small neighborhood in $\overline\Omega$.
For convenience, we still use the notation $\hat u$ for the extension.
For $x\notin\Omega$, $x_n\ge0$, let $z=(x',\phi(x'))$.
Then we have
\[
\begin{aligned}
\hat u(x)&\ge\hat u(z)+D\hat u(z)\cdot(x-z)\ge c|z|^q-C|z|^{1+\beta+a}\ge c|x|^q,\\
\hat u(x)&\le\hat u(z)\le C|x|^{1+a}.
\end{aligned}
\]
Here $1+\beta+a>q$.
Together with \eqref{initialgrowth}, we obtain
\[
c|x|^q\le\hat u(x)\le C|x|^{1+a}\qquad(x_n\ge0).
\]
It follows by difference quotients that
\begin{equation}\label{extensiongradient}
|p|\le C r^a\qquad\bigl(p\in\partial\hat u(x),\ x_n\ge0,\ |x|\le r\bigr).
\end{equation}

Define
\[
D_h'=\operatorname{proj}_{e_n^\perp}\{\hat u<h,\ x_n>0\},\qquad a_h=\max_{\overline{D_h'}}\phi,
\]
\[
D_h^+=\{\hat u<h,\ x_n>a_h\},\qquad D_h=\{x\mid \hat u(x',a_h+|x_n-a_h|)<h\}.
\]
The upper part $D_h^+$ lies in $\Omega$.
The function $x\mapsto\hat u(x',a_h+|x_n-a_h|)$ is convex since $\hat u$ is convex and nondecreasing in $x_n$.
Hence $D_h$ is convex and symmetric about $x_n=a_h$.
The growth estimates imply
\[
R_h=C h^{b_-},\qquad r_h=c h^{b_+},\qquad a_h\le C R_h^{1+\beta}=o(r_h),\qquad B_{c r_h}\subset D_h\subset B_{C R_h}.
\]

For $x=(x',a_h)\in\overline{D_h^+}$ and $z=(x',\phi(x'))$,
\[
T_n(z)=\phi^*(T'(z))\le C|T'(z)|^{1+\beta},\qquad |T(x)-T(z)|\le C a_h^a.
\]
It follows that
\[
0\le\hat u_n(x)\le C R_h^{a(1+\beta)}=j_h.
\]
Indeed we have
\[
\hat u_n\le C\bigl(R_h^{a(1+\beta)}+a_h^a+R_h^{1+2a}\bigr)\le C R_h^{a(1+\beta)}.
\]

On $D_h^+$ the equation is
\[
\det(D^2\hat u-E)=\rho,\qquad E(x)=c_{xx}(x,T(x))+D^2\varphi(x),\qquad \rho(x)=|\det c_{xy}(x,T(x))|\frac{f(x)}{g(T(x))}.
\]
Write $\rho_0=f(0)/g(0)$.
Then we obtain
\[
\|E\|_{L^\infty(D_h^+)}\le C R_h^{2a}=e_h,\qquad \|\rho/\rho_0-1\|_{L^\infty(D_h^+)}\le C R_h^{a\alpha}=d_h.
\]
Let $w_h$ be the Aleksandrov solution of
\[
\det D^2w_h=\rho_0\quad\hbox{in }D_h,\qquad w_h=h\quad\hbox{on }\partial D_h.
\]
Existence and uniqueness follow from \cite[Theorem 2.1]{TW2008}.
Its interior regularity follows from \cite[Lemma 3.4 and Theorem 3.1]{TW2008}.
Let $H_h=h-\min_{D_h}w_h$.
It follows from symmetry that
\[
(w_h)_n=0\qquad\hbox{on }D_h\cap\{x_n=a_h\}.
\]

Take the fixed constant in $R_h$ large so that $D_h\subset B_{R_h}(a_he_n)$ and define
\[
Q_h(x)=|x-a_he_n|^2-R_h^2,\qquad s_\pm=(1\pm2d_h)^{1/n}.
\]
For $\eta>0$, we define two test functions
\[
\begin{aligned}
w_-&=s_+(w_h-h)+h+e_hQ_h+(j_h+\eta)(x_n-a_h-R_h),\\
w_+&=s_-(w_h-h)+h-e_hQ_h-\eta(x_n-a_h-R_h).
\end{aligned}
\]
On the curved part of $\partial D_h^+$, we have
\[
w_-\le h=\hat u\le w_+.
\]
On the flat part, we have
\[
(w_-)_n=j_h+\eta>\hat u_n\ge0>(w_+)_n=-\eta.
\]
A positive maximum of $w_--\hat u$ cannot be reached on the flat part.
At an interior maximum, we obtain that
\[
D^2\hat u-E\ge D^2w_--E\ge s_+D^2w_h+e_hI>0,
\]
\[
\rho=\det(D^2\hat u-E)\ge\det(s_+D^2w_h+e_hI)\ge(1+2d_h)\rho_0>\rho,
\]
which is impossible.
Similarly, at an interior positive maximum of $\hat u-w_+$,
\[
0<D^2\hat u-E\le D^2w_+-E\le s_-D^2w_h-e_hI,
\]
\[
\rho\le\det(s_-D^2w_h-e_hI)\le(1-2d_h)\rho_0<\rho,
\]
which is also impossible.
Thus we have $w_-\le\hat u\le w_+$ on $D_h^+$.
Letting $\eta\to0$ implies
\begin{equation}\label{matrixcomparison}
\|\hat u-w_h\|_{L^\infty(D_h^+)}\le C\bigl(d_hH_h+e_hR_h^2+j_hR_h\bigr).
\end{equation}

It follows from the choice of $\varepsilon$ that
\[
\delta>2\varepsilon,\qquad 2a-1>\delta,\qquad b_-(1+\beta+a)-1>\delta,
\]
\[
d_h\le C h^\delta,\qquad e_hR_h^2+j_hR_h\le C h^{1+\delta},\qquad \hat u(a_he_n)\le C a_h^{1+a}=o(h).
\]
A minimum of $w_h$ lies on its symmetry plane, where $\hat u\ge0$.
By \eqref{matrixcomparison},
\[
H_h\le h+C d_hH_h+o(h),\qquad H_h\ge h-Cd_hH_h-o(h).
\]
Therefore we have
\[
c h\le H_h\le C h.
\]
The volume estimate \cite[Lemma 3.2]{TW2008} shows that
\[
c H_h^{n/2}\le |D_h|\le C H_h^{n/2}.
\]
In particular,
\[
c h^{n/2}\le |D_h|\le C h^{n/2},\qquad \|\hat u-w_h\|_{L^\infty(D_h^+)}\le C h^{1+\delta}.
\]

Let $z_h(x)=(x',a_h+|x_n-a_h|)$.
For $x\in D_{h/4}$, we have
\[
z_h(x)-z_{h/4}(x)=t e_n,\qquad 0\le t\le2a_h.
\]
Note that the segment joining these points lies in $\overline\Omega$. It follows by \eqref{extensiongradient} that
\[
|\hat u(z_h(x))-\hat u(z_{h/4}(x))|\le C R_h^aa_h\le C R_h^{1+\beta+a}\le C h^{1+\delta}.
\]
Hence $D_{h/4}\subset D_h$ for sufficiently small $h$.
The reflection symmetry and the two comparison inequalities give
\[
\|w_h-w_{h/4}\|_{L^\infty(D_{h/4})}\le C h^{1+\delta}.
\]
If $x\in D_h\cap\Omega$ lies below the symmetry plane, the same estimate on the segment $[x,z_h(x)]$ compares $\hat u(x)$ with $w_h(x)$.
This proves \eqref{modelerror}.
Finally, we have
\begin{equation}\label{modelorigin}
|w_h(0)|=|w_h(2a_he_n)|\le C h^{1+\delta}+C a_h^{1+a}\le C h^{1+\delta}.\qedhere
\end{equation}
\end{proof}

\begin{proposition}\label{boundaryjets}
At every boundary contact pair, in the coordinates \eqref{cnorm}, there is a positive definite matrix $H_b$ such that
\begin{equation}\label{boundaryexpansion}
\left|u(x)-\frac12x^TH_bx\right|\le C|x|^{2+\gamma},\qquad cI\le H_b\le CI,\qquad \det H_b=\frac{f(0)}{g(0)}.
\end{equation}
The analogous estimate holds for $v$.
\end{proposition}

\begin{proof}
Let $\widehat L_h$ be a John matrix for the reflected domain $D_h$.
It follows from Lemma~\ref{boundarycomparison} that
\[
\|\widehat L_h^{-1}\|\le C h^{-b_+},\qquad \|\widehat L_h^{-1}\widehat L_{h/4}\|+\|\widehat L_{h/4}^{-1}\widehat L_h\|\le C.
\]
By \eqref{modelorigin}, the origin lies in a fixed negative sublevel of the normalized functions $(w_h-h)/h$ and $(w_{h/4}-h/4)/(h/4)$.
The Aleksandrov maximum estimate \cite[Lemma 3.3]{TW2008} shows that it has a fixed distance away from the normalized boundaries.
Hence by the constant Monge--Amp\`ere interior estimates \cite[Lemma 3.6 and Theorem 3.1]{TW2008} we obtain uniform positive Hessian bounds and $C^4$ bounds on a common normalized ball about the origin.

The difference $w_h-w_{h/4}$ satisfies
\[
\left(\int_0^1\bigl[(1-t)D^2w_h+tD^2w_{h/4}\bigr]^{-1}_{ij}\,dt\right)D_{ij}(w_h-w_{h/4})=0.
\]
In the normalized coordinates this equation is uniformly elliptic and its coefficients have uniformly bounded $C^2$ norms.
By differentiating it, we use the Schauder estimate \cite[Theorem 6.2]{GT2001} and \eqref{modelerror} to obtain that
\begin{equation}\label{modelderivatives}
\sup_{B_{c h^{b_+}}}|D^j(w_h-w_{h/4})|\le C h^{1+\delta-jb_+}\qquad(j=0,1,2,3).
\end{equation}
For $h_k=4^{-k}h_0$,
\[
\sum_{k=0}^\infty |D^2w_{h_k}(0)-D^2w_{h_{k+1}}(0)|\le C\sum_{k=0}^\infty h_k^{\delta-2\varepsilon}<\infty.
\]
Consequently,
\[
D^2w_{h_k}(0)\longrightarrow H_b,\qquad \det H_b=\rho_0,\qquad cI\le H_b\le CI.
\]
It follows by convexity that
\[
|Dw_h(0)|\le C h^{1-b_+}\longrightarrow0.
\]
Summing \eqref{modelderivatives} yields
\[
|D^2w_h(0)-H_b|\le C h^{\delta-2\varepsilon},\qquad |Dw_h(0)|\le C h^{1+\delta-b_+},
\]
\[
\sup_{B_{c h^{b_+}}}|D^3w_h|\le C h^{1+\delta-3b_+}.
\]
Here and below it is enough to take $h=h_k$.
By Taylor's formula, \eqref{modelerror} and \eqref{modelorigin}, we have
\[
\left|\hat u(x)-\frac12x^TH_bx\right|\le C h^{1+\delta}\qquad\bigl(x\in\overline\Omega\cap B_{c h^{b_+}}\bigr).
\]
Since $\varphi(x)=O(|x|^{2+2a})$ and $b_+(2+2a)=2$, we obtain
\[
\left|u(x)-\frac12x^TH_bx\right|\le C|x|^{2+\gamma_0},\qquad \gamma_0=\frac{\delta-2\varepsilon}{b_+}>0.
\]
The constants are uniform in the boundary pair.
Applying the same argument to $v$ gives the analogous estimate.

For $y=T(x)$, the two expansion estimate and $c(x,y)=x\cdot y+O(|x|^2|y|^2)$ imply
\[
c\bigl(|x|^2+|y|^2\bigr)\le u(x)+v(y)=c(x,y)\le C|x||y|.
\]
Thus, in the original coordinates, we have
\begin{equation}\label{boundarylipschitz}
|T(x)-T(b)|\le C|x-b|\qquad\bigl(b\in\partial\Omega,\ x\in\overline\Omega\bigr).
\end{equation}
The inverse map satisfies the analogous estimate.
The corrected potential and its extension now satisfy
\[
c|x|^2\le\hat u(x)\le C|x|^2,\qquad |p|\le Cr\quad\bigl(p\in\partial\hat u(x),\ x_n\ge0,\ |x|\le r\bigr).
\]
Indeed, below $\partial\Omega$ the change from $x$ to $(x',\phi(x'))$ contributes at most $C|x|^{2+\beta}$.

Repeat the construction in Lemma~\ref{boundarycomparison}, with the same function $\varphi$ in \eqref{correction}.
Writing $r=\sqrt h$, we have
\[
B_{c r}\subset D_h\subset B_{C r},\qquad a_h\le C r^{1+\beta}.
\]
It follows from the boundary estimate \eqref{boundarylipschitz} that, on the cut,
\[
T_n(x)\le C r^{1+\beta},\qquad 0\le\hat u_n(x)\le C\bigl(r^{1+\beta}+r^{1+2a}\bigr).
\]
The quantities in \eqref{matrixcomparison} therefore satisfy
\[
j_h\le C\bigl(r^{1+\beta}+r^{1+2a}\bigr),\qquad d_h\le C r^\alpha,\qquad e_h\le C r^{2a}.
\]
Since $2a>1>\gamma$,
\[
d_hh+e_hr^2+j_hr\le C\bigl(r^{2+\alpha}+r^{2+\beta}+r^{2+2a}\bigr)\le C r^{2+\gamma}.
\]
The change under reflection is at most
\[
C r a_h\le C r^{2+\beta}\le C r^{2+\gamma}.
\]
Consequently, we obtain
\[
\|\hat u-w_h\|_{L^\infty(D_h\cap\Omega)}+\|w_h-w_{h/4}\|_{L^\infty(D_{h/4})}+|w_h(0)|\le C h^{1+\gamma/2}.
\]
The same derivative estimates are now
\[
\sup_{B_{c\sqrt h}}|D^j(w_h-w_{h/4})|\le C h^{1+\gamma/2-j/2}\qquad(j=0,1,2,3).
\]
By summation we obtain
\[
|D^2w_h(0)-H_b|\le C h^{\gamma/2},\qquad |Dw_h(0)|\le C h^{(1+\gamma)/2},\qquad \sup_{B_{c\sqrt h}}|D^3w_h|\le C h^{(\gamma-1)/2}.
\]
Finally, Taylor's formula proves \eqref{boundaryexpansion}.
\end{proof}

\begin{proof}[Proof of Theorem~\ref{HR}]
Restore the smooth coordinate changes and the supporting functions in Proposition~\ref{boundaryjets}.
At each $b\in\partial\Omega$ there is an ordinary quadratic polynomial $P_b$ with
\begin{equation}\label{ordinaryjets}
|u(x)-P_b(x)|\le C|x-b|^{2+\gamma},\qquad P_b(b)=u(b),\qquad DP_b(b)=Du(b).
\end{equation}
If $b,b'\in\partial\Omega$ and $r=|b-b'|$ is small, the $C^1$ boundary admits a common interior ball of radius $c r$ within distance $C r$ of both points.
On this ball,
\[
|P_b-P_{b'}|\le C r^{2+\gamma}.
\]
The elementary norm equivalence for quadratic polynomials gives
\begin{equation}\label{jetcompatibility}
|D^2P_b-D^2P_{b'}|\le C|b-b'|^\gamma.
\end{equation}

Let $z\in\Omega$ and $d=\operatorname{dist}(z,\partial\Omega)$.
Choose a nearest boundary point $b$ and return to the fixed coordinates \eqref{cnorm} at $(b,T(b))$.
The point $b$ is now the origin and $P_b(x)=x^TH_bx/2$.
The coordinate change gives
\[
|z|\le C d,\qquad B_{c d}(z)\subset\Omega,\qquad \|u-P_b\|_{L^\infty(B_{c d}(z))}\le C d^{2+\gamma}.
\]
Corollary~\ref{interiorholder} gives $u\in C^{2,\alpha}_{\rm loc}(\Omega)$.
We next make its interior estimate uniform as $d\to0$.

Put $H=H_b$, choose $r=c_0d$ with fixed sufficiently small $c_0>0$, and set
\[
x=z+rH^{-1/2}X,\qquad y=Hz+rH^{1/2}Y,
\]
\[
\tilde u(X)=r^{-2}\bigl(u(x)-P_b(z)-Hz\cdot(x-z)\bigr),
\]
\[
\tilde c(X,Y)=r^{-2}\bigl(c(x,y)-c(z,y)-Hz\cdot(x-z)\bigr).
\]
The source coordinate contains $B_4$.
Dividing both marginals by $r^n(\det H)^{-1/2}f(0)$, we have
\[
\tilde f(X)=\frac{f(x)}{f(0)},\qquad \tilde g(Y)=\frac{g(y)}{g(0)},
\]
because $\det H=f(0)/g(0)$.
The focused cost estimates and \eqref{boundaryexpansion} give, on fixed balls,
\[
\left\|\tilde u-\frac12|X|^2\right\|_{C^0}\le C d^\gamma,\qquad \|\tilde c-X\cdot Y\|_{C^{2,\alpha}}\le C d^2,
\]
\[
\|\tilde f-1\|_{C^0}+\|\tilde g-1\|_{C^0}\le C d^\alpha,\qquad \|\tilde f\|_{C^\alpha}+\|\tilde g\|_{C^\alpha}\le C.
\]
By \eqref{boundarylipschitz}, all actual normalized targets lie in a fixed ball.
The semiconvexity estimate and the first displayed inequality imply
\[
\|\tilde T-\operatorname{id}\|_{L^\infty(B_2)}\le C\bigl(d^{\gamma/2}+d^2\bigr).
\]
It follows by Proposition~\ref{contacts} that
\[
C_1=\overline B_1,\qquad C_2=\tilde T(C_1),\qquad B_{1/3}\subset C_1,C_2\subset B_3.
\]
The restriction of the full transport plan to $C_1\times C_2$ has precisely the restricted source and target marginals.
This is the localization in \cite[Section 3]{DPF2015}.
The perturbation estimate \cite[Theorem 5.3]{DPF2015} yields, on a fixed smaller ball,
\begin{equation}\label{interioruniform}
\|\tilde u\|_{C^{2,\alpha}}\le C,\qquad cI\le D^2\tilde u-\tilde c_{XX}(X,\tilde T(X))\le CI.
\end{equation}
The constants are independent of $d$.

On a fixed small ball in the variables $X=(x-z)/d$, set
\[
\tilde w(X)=d^{-2}\bigl(u(z+dX)-P_b(z+dX)\bigr),
\]
\[
A_d(X)=c_{xx}(z+dX,T(z+dX)),\qquad \rho_d(X)=\rho(z+dX).
\]
By \eqref{ordinaryjets},
\[
\|\tilde w\|_{C^0}\le C d^\gamma.
\]
In the focused coordinates,
\[
|c_{xx}(x,y)|\le C|y|^2,\qquad |c_{xxx}(x,y)|+|c_{xxy}(x,y)|\le C(|x|+|y|).
\]
Together with \eqref{interioruniform} and $|T(x)|\le C d$, these give
\[
\|A_d\|_{C^\gamma}\le C d^2,\qquad \|\log\rho_d-\log\rho_0\|_{C^\gamma}\le C d^\alpha.
\]
Define
\[
a^{ij}(X)=\int_0^1\left[(1-t)H_b+t\bigl(D^2u(z+dX)-A_d(X)\bigr)\right]^{-1}_{ij}\,dt.
\]
By the equation and \eqref{interioruniform} we have
\[
a^{ij}\tilde w_{ij}=\log\rho_d-\log\rho_0+a^{ij}(A_d)_{ij},
\]
\[
cI\le(a^{ij})\le CI,\qquad \|(a^{ij})\|_{C^\gamma}\le C.
\]
The interior Schauder estimate \cite[Theorem 6.2]{GT2001} now shows
\[
\|\tilde w\|_{C^{2,\gamma}(B_\eta)}\le C d^\gamma.
\]
Returning to the original scale, we obtain
\begin{equation}\label{boundaryhessian}
|D^2u(z)-D^2P_b|\le C d^\gamma,\qquad [D^2u]_{C^\gamma(B_{\eta d}(z))}\le C.
\end{equation}
The smooth coordinate changes preserve these estimates.

For two points at distance at most a fixed small multiple of their distance to the boundary, use the second estimate in \eqref{boundaryhessian}.
Otherwise, use the first estimate at their nearest boundary points and \eqref{jetcompatibility}.
Hence
\[
u\in C^{2,\gamma}(\overline\Omega),\qquad \|u\|_{C^{2,\gamma}(\overline\Omega)}\le C.
\]
Applying the same argument to the dual potential, we obtain
\[
v\in C^{2,\gamma}(\overline{\Omega^*}).
\]
This proves Theorem~\ref{HR}.
\end{proof}

\section*{Acknowledgments}

We began studying this problem when the first author joined Westlake University as a postdoctoral researcher. Inspired by Collins and Tong's July 2025 arXiv preprint, we made substantial progress by establishing a monotonicity inequality for general costs by the end of 2025. In late July 2026, we completed the proof under A3s. Believing the method might also apply under A3w, we postponed finalizing the A3s manuscript to pursue the A3w condition.

In August 2026, the first author accepted an invitation to speak at the Feishu Mathematics Symposium in Shanghai. On September 4, the first author emailed the organizers the title, \emph{Global $C^{2,\alpha}$ regularity for optimal transport under A3s condition}, and the abstract. In the September 29 talk, we also mentioned our ongoing work under A3w.

With AI assistance, we completed the proof of version 1 under A3w assuming the domains are mutually uniformly $c$-convex with $C^2$ boundaries. Through continued discussions with GPT-6 Astra, we have since substantially reworked the proof, developing alternative arguments and simplifying its overall structure, and extending the results to mutually $c$-convex domains.

We used ChatGPT (OpenAI) in preparing this manuscript. The authors have reviewed all AI-assisted content and take full responsibility for the final manuscript.

\end{document}